\documentclass[11pt]{article}
\usepackage[a4paper,margin=1in]{geometry}
\usepackage{amsmath,amssymb,amsthm,bm}
\usepackage{graphicx}
\usepackage{hyperref}
\usepackage{physics}
\usepackage{cite}
\usepackage{enumitem}
\usepackage{booktabs,siunitx}
\usepackage{subcaption}

\AtBeginDocument{\RenewCommandCopy\qty\SI}

\newcommand{\PP}{\mathbb{P}}

\newcommand{\PII}{\mathrm{P_{II}}}

\newcommand{\PIIIp}{\mathrm{P_{III}'}}
\newcommand{\PIV}{\mathrm{P_{IV}}}
\newcommand{\PVI}{\mathrm{P_{VI}}}

\DeclareGraphicsExtensions{.pdf,.png,.jpg,.jpeg}
\newcommand{\imgw}[2][]{\includegraphics[width=\linewidth,#1]{#2}}

\newcommand{\Kop}{\mathsf{K}}   
\newcommand{\Ai}{\mathrm{Ai}}   

\usepackage[T1]{fontenc}
\usepackage[utf8]{inputenc}

\newtheorem{theorem}{Theorem}

\newtheorem{remark}{Remark}

\DeclareUnicodeCharacter{03C3}{\ensuremath{\sigma}}   
\DeclareUnicodeCharacter{2032}{\ensuremath{^\prime}}  
\DeclareUnicodeCharacter{2013}{--}                    

\newcommand{\R}{\mathbb{R}}
\newcommand{\C}{\mathbb{C}}

\title{\textbf{The Instability of Painlevé Equations in Recovering Largest Eigenvalue Distributions of GUE, LUE, JUE and an Attempt of Solution to It}}
\author{Haonan Gu}
\date{November 2025}

\begin{document}
\maketitle

\begin{abstract}
The distribution of the largest eigenvalue for the three classical unitary ensembles---GUE, LUE, and JUE---admits two complementary exact descriptions: (i) as Fredholm determinants of their orthogonal--polynomial correlation kernels and (ii) as isomonodromic $\tau$--functions governed by Painlevé equations. For finite $n$, the associated Jimbo--Miwa--Okamoto $\sigma$--forms are $\PIV$ (GUE), $\mathrm{PV}$ (LUE), and $\PVI$ (JUE); under soft- or hard-edge scalings these degenerate to $\PII$ or $\PIIIp$ descriptions of the Tracy--Widom and hard-edge laws \cite{tracy1994level,forrester2003painleve,deift1999orthogonal}. 

It is well known among random matrix theorists (for example Folkmar Bornemann) that the Fredholm determinant is a more numerically stable and accurate way to compute the CDF of the largest eigenvalue for GUE, LUE, JUE than direct Painlevé integration. The aim of this paper is not to improve on Fredholm methods, but to see to what extent one can numerically recover the \emph{correct} Painlevé solution from finite-$n$ data and how unstable this reconstruction is. Numerically, we verify the equality between the Fredholm- and Painlevé-based CDFs by combining (a) high-accuracy Nyström discretizations of the finite-$n$ Fredholm determinants \cite{bornemann2010numerical} with (b) an anchored, branch-locked integration of the $\sigma$--form ODEs, where anchors are extracted from local least-squares fits to $\log\det(I-\mathsf K)$. Our results confirm agreement across GUE/LUE/JUE with precision of $O(10^{-3})$ to $O(10^{-5})$ (occasionally $O(10^{-2})$) and illustrate the finite-$n$ to scaling-limit transition. The theoretical connections to $\tau$--functions and Virasoro constraints follow the framework of \cite{adler2000random,forrester2003painleve}.
\end{abstract}

\section{Introduction}

\subsection{Preliminaries: ensembles, kernels, Fredholm determinants, Painlevé connections}

We work with the classical unitary ensembles associated with Hermite, Laguerre, and Jacobi orthogonal polynomials \cite{deift1999orthogonal}. In the Gaussian unitary ensemble (GUE) the eigenvalues are supported on $\R$ with weight $w_{\mathrm H}(x)=e^{-x^2}$. In the Laguerre unitary ensemble (LUE) the support is $[0,\infty)$ and the weight is $w_{\mathrm L}(x)=x^{\alpha} e^{-x}$ with $\alpha>-1$. In the Jacobi unitary ensemble (JUE) the support is $(-1,1)$ and the weight is $w_{\mathrm J}(x)=(1-x)^{b}(1+x)^{a}$ with $a,b>-1$.

For each ensemble we write $\{\phi_k\}_{k\ge0}$ for the orthonormal polynomial functions
\[
\phi_k(x)=p_k(x)\sqrt{w(x)},\qquad 
\int \phi_k(x)\phi_\ell(x)\,dx=\delta_{k\ell},
\]
where $p_k$ is the orthonormal polynomial (Hermite, Laguerre, Jacobi) for the corresponding weight $w$.

The eigenvalues form a determinantal point process with $m$--point correlations
\[
\rho_m(x_1,\dots,x_m)=\det\!\big[K^{(n)}(x_i,x_j)\big]_{i,j=1}^m,
\]
where the finite-$n$ correlation kernel is the orthogonal projection
\[
K^{(n)}(x,y)=\sum_{k=0}^{n-1}\phi_k(x)\phi_k(y).
\]
For numerical work we use the Christoffel--Darboux (CD) representation together with a diagonal fallback:
\[
K^{(n)}(x,y)=
\frac{\gamma_n\big(\phi_{n-1}(x)\phi_n(y)-\phi_n(x)\phi_{n-1}(y)\big)}{x-y},
\qquad 
K^{(n)}(x,x)=\sum_{k=0}^{n-1}\phi_k(x)^2,
\]
with ensemble-dependent constants $\gamma_n$; see \cite[Ch.~2--5]{deift1999orthogonal}.

If $B$ is a Borel set in the support and
\[
(\Kop f)(x)=\int_B K^{(n)}(x,y)f(y)\,dy\qquad\text{on }L^2(B),
\]
then the probability that $B$ contains no eigenvalue is a Fredholm determinant
\[
E(B)=\det(I-\Kop)_{L^2(B)}.
\]
Taking $B$ to be the one-sided interval at the upper edge gives the largest-eigenvalue CDFs:
\[
\text{GUE: }F_n(s)=\det(I-\Kop)_{(s,\infty)},\quad
\text{LUE: }F_{n,\alpha}(s)=\det(I-\Kop)_{(s,\infty)},\quad
\text{JUE: }F_{n,a,b}(s)=\det(I-\Kop)_{(s,1)}.
\]
We evaluate these determinants numerically via Nyström discretization \cite{bornemann2010numerical}. For a quadrature rule $\{t_i,w_i\}_{i=1}^N$ on $B$ we form the symmetric matrix
\[
A_{ij}=\sqrt{w_i}\,K^{(n)}(t_i,t_j)\sqrt{w_j},\qquad
F=\det(I-A)=\exp\!\Big(\sum_{j=1}^N \log \lambda_j\Big),
\]
where $\{\lambda_j\}$ are the eigenvalues of $I-A$.

These Fredholm determinants are also isomonodromic $\tau$--functions and satisfy Painlevé equations with ensemble-specific parameters \cite{adler2000random,forrester2003painleve}. For finite $n$ the largest-eigenvalue gaps are described by a $\sigma$–form of $\PIV$ for GUE, a $\sigma$–form of $\mathrm{PV}$ for LUE, and a $\sigma$–form of $\PVI$ for JUE. At the soft edge these reduce to $\PII$ (Tracy–Widom), while hard-edge scalings produce $\PIIIp$ (Bessel kernel) \cite{tracy1994level,forrester2003painleve,deift1999orthogonal}. One central theme of this paper is that, although these ODE characterizations are exact, the corresponding numerical IVPs are highly unstable and must be locked to the Fredholm side.

\subsection{Finite-$n$ \texorpdfstring{$\PIV$}{PIV} for the GUE gap}\label{sec:finite-n-PIV}

Let $F_n$ denote the finite-$n$ GUE largest-eigenvalue distribution
\[
F_n(s) := \PP(\lambda_{\max}\le s),
\]
for the $n\times n$ GUE with weight $e^{-x^2}$. Equivalently,
\[
F_n(s)=\PP(\lambda_{\max}\le s)
      = \PP(\text{no eigenvalues lie in }(s,\infty)).
\]

Following Forrester and Witte \cite{forrester2003painleve}, write $E_n(0;(s,\infty))$ for the probability that $(s,\infty)$ is free of eigenvalues. Their analysis (see Eq.~(1.4)--(1.7) of \cite{forrester2003painleve}) shows that there is a distinguished function $R_n(s)$ such that
\begin{equation}
F_n(s)=E_n(0;(s,\infty))
= \exp\!\Big(-\int_s^\infty R_n(t)\,dt\Big),
\label{eq:Fn-Rn-int}
\end{equation}
and $R_n$ satisfies the Jimbo--Miwa--Okamoto $\sigma$--form of $\PIV$:
\begin{equation}
\big(R_n''(s)\big)^2
=4\big(s\,R_n'(s)-R_n(s)\big)^2
-4\big(R_n'(s)\big)^2\big(R_n'(s)+2n\big).
\label{eq:PIV-sigma-R}
\end{equation}
The large-$s$ behaviour is fixed by the Hermite-kernel tail:
\begin{equation}
R_n(s)\sim
\frac{2^{\,n-1}}{\sqrt{\pi}(n-1)!}\,s^{2n-2}e^{-s^2},
\qquad s\to+\infty.
\label{eq:Rn-asymp}
\end{equation}

Define the logarithmic derivative
\begin{equation}
\sigma_n(s):=\frac{d}{ds}\log F_n(s).
\label{eq:sigma-def}
\end{equation}
Differentiating \eqref{eq:Fn-Rn-int} and using
$\frac{d}{ds}\int_s^\infty R_n(t)\,dt=-R_n(s)$ gives
\begin{equation}
\sigma_n(s)=\frac{d}{ds}\log F_n(s)=R_n(s).
\label{eq:sigma-Rn-ident}
\end{equation}
Thus the Forrester--Witte function $R_n$ is exactly the logarithmic derivative of $F_n$.

Combining \eqref{eq:PIV-sigma-R} and \eqref{eq:sigma-Rn-ident} yields the $\sigma$--form $\PIV$ equation for $\sigma_n$.

\begin{theorem}[Jimbo--Miwa--Okamoto $\sigma$--form for the finite-$n$ GUE]
\label{thm:Fn-sigma-PIV}
Let $F_n(s)=\PP(\lambda_{\max}\le s)$ be the largest-eigenvalue CDF of the $n\times n$ GUE with weight $e^{-x^2}$, and define $\sigma_n$ by \eqref{eq:sigma-def}. Then:
\begin{enumerate}[label=\emph{(\roman*)},leftmargin=2em]
  \item $\sigma_n$ satisfies the Jimbo--Miwa--Okamoto $\sigma$--form
        \begin{equation}
        \big(\sigma_n''(s)\big)^2
        =4\big(s\,\sigma_n'(s)-\sigma_n(s)\big)^2
        -4\big(\sigma_n'(s)\big)^2\big(\sigma_n'(s)+2n\big).
        \label{eq:PIV-sigma}
        \end{equation}
  \item As $s\to+\infty$ one has
        \begin{equation}
        \sigma_n(s)\sim
        \frac{2^{\,n-1}}{\sqrt{\pi}(n-1)!}\,s^{2n-2}e^{-s^2},
        \label{eq:sigma-asymp}
        \end{equation}
        so in particular $\sigma_n(s)\to0$.
  \item Conversely, any smooth solution $\sigma_n$ of
        \eqref{eq:PIV-sigma} with asymptotics \eqref{eq:sigma-asymp} reconstructs the CDF via
        \begin{equation}
        F_n(s)=\exp\!\Big(\int_{s_0}^{s}\sigma_n(u)\,du\Big),
        \qquad
        \lim_{s\to+\infty}F_n(s)=1,
        \label{eq:PIV-to-Fn}
        \end{equation}
        for any reference point $s_0$.
\end{enumerate}
\end{theorem}

\begin{remark}
Equation \eqref{eq:PIV-sigma} is the special case $a=0$ of the general $\sigma$--$\PIV$ satisfied by the spectral averages $\tilde{E}_N(\lambda;a)$ in Forrester--Witte \cite[Eq.~(1.9), Eq.~(4.15)]{forrester2003painleve}, after identifying $\sigma_n(s)$ with $U_N(s;0)$ in their notation and setting $N=n$.
\end{remark}

For completeness we recall the relation between the $\sigma$--form \eqref{eq:PIV-sigma} and the usual $\PIV$
\begin{equation}
y''=\frac{(y')^2}{2y}+\frac{3}{2}y^3+4xy^2+2(x^2-\alpha)y+\frac{\beta}{y},
\label{eq:PIV-standard}
\end{equation}
in Okamoto's Hamiltonian formulation. Introduce canonical coordinates $q(t),p(t)$ and the Hamiltonian
\begin{equation}
H(t;q,p)
 = (2p-q-2t)\,p\,q - 2\alpha_1 p - \alpha_2 q,
\label{eq:Okamoto-H}
\end{equation}
where $(\alpha_0,\alpha_1,\alpha_2)$ satisfy $\alpha_0+\alpha_1+\alpha_2=1$. Hamilton's equations are
\begin{equation}
\begin{aligned}
q' &= q(4p - q - 2t) - 2\alpha_1,\\
p' &= 2p(q + t - p) + \alpha_2.
\end{aligned}
\label{eq:qp-Okamoto}
\end{equation}
Eliminating $p$ yields \eqref{eq:PIV-standard} with parameters $(\alpha,\beta)$ linear in $(\alpha_0,\alpha_1,\alpha_2)$. For the GUE spectral averages $\tilde{E}_N(\lambda;a)$ one finds \cite[Eq.~(4.10), Eq.~(4.15), Eq.~(4.46)--(4.47)]{forrester2003painleve}
\[
\alpha_0 = 1+N+a,\qquad \alpha_1 = -a,\qquad \alpha_2=-N.
\]
In our largest-eigenvalue problem ($a=0$, $N=n$) this becomes
\[
\alpha_0 = 1+n,\qquad \alpha_1 = 0,\qquad \alpha_2=-n,
\]
and \eqref{eq:PIV-standard} has $(\alpha,\beta)=(2n+1,0)$. The logarithmic derivative $\sigma_n$ is an affine transform of the Hamiltonian:
\begin{equation}
\sigma_n(s) = U_n(s;0)
            = H\bigl(s;q_n(s),p_n(s)\bigr) + \text{(explicit affine term)},
\label{eq:sigma-H}
\end{equation}
where $(q_n,p_n)$ solve \eqref{eq:qp-Okamoto}. The precise affine term, recorded in \cite{forrester2003painleve}, plays no role numerically because it disappears from the $\sigma$--form.

\subsection{Edge limit: from \texorpdfstring{$\PIV$}{PIV} to \texorpdfstring{$\PII$}{PII} (Tracy--Widom)}

At the upper edge $\sqrt{2n}$ set
\[
s=\sqrt{2n}+x\,n^{1/6},\qquad t=\sqrt{2n}+y\,n^{1/6}.
\]
Then, uniformly for $(x,y)$ in compact sets,
\[
n^{1/6}K^{(n)}_{\mathrm H}(s,t)\ \longrightarrow\ K_{\mathrm{Ai}}(x,y)
=\frac{\Ai(x)\Ai'(y)-\Ai'(x)\Ai(y)}{x-y}.
\]
Consequently, the soft-edge limit of the largest-eigenvalue CDF exists:
\begin{equation}
F_2(x):=\lim_{n\to\infty}F_n\!\big(\sqrt{2n}+x\,n^{-1/6}\big)
=\det\nolimits_{L^2([x,\infty))}\!\big(I-\Kop_{\mathrm{Ai}}\big).
\label{eq:Fredholm}
\end{equation}

Let $q$ be the Hastings--McLeod solution of $\PII$
\begin{equation}
q''(x)=x\,q(x)+2q(x)^3,
\qquad
q(x)\sim \Ai(x)\quad(x\to+\infty).
\label{eq:PII}
\end{equation}
This solution is real, pole-free on $\R$, and behaves like $q(x)\sim \sqrt{-x/2}$ as $x\to-\infty$. Defining
\begin{equation}
F_2(x)
=\exp\!\Big[-\!\int_x^\infty (t-x)\,q(t)^2\,dt\Big],
\label{eq:TW-PII}
\end{equation}
Tracy and Widom showed \cite{tracy1994level} that the function $F_2$ in \eqref{eq:TW-PII} coincides with the Fredholm determinant in \eqref{eq:Fredholm}. Thus $F_2$ is the limiting distribution of the scaled GUE largest eigenvalue and admits both a Fredholm and a Painlevé representation.

For finite $n$, the map $F_n\mapsto \sigma_n=\dv{}{s}\log F_n$ solves \eqref{eq:PIV-sigma}, and integrating \eqref{eq:PIV-to-Fn} recovers $F_n$. Under soft-edge scaling, $F_n$ converges to $F_2$, which in turn satisfies both the Airy-kernel Fredholm representation \eqref{eq:Fredholm} and the Painlevé--II representation \eqref{eq:TW-PII}. One can therefore track the finite-$n$ GUE edge all the way to the universal Tracy--Widom regime through a chain of Painlevé equations ($\PIV\to\PII$) and kernels (Hermite $\to$ Airy).

\subsection{Finite-\texorpdfstring{$N$}{N} LUE: Fredholm determinants and $\sigma$--\texorpdfstring{$\mathrm{PV}$}{PV}}\label{sec:LUE-PV}

For $a>-1$, the LUE consists of $N\times N$ Hermitian positive-definite matrices with joint eigenvalue density
\[
\frac{1}{Z_{N,a}} \prod_{1\le i<j\le N}|\lambda_i-\lambda_j|^2 \prod_{i=1}^N \lambda_i^{a}e^{-\lambda_i}\mathbf{1}_{(0,\infty)}(\lambda_i)\,d\lambda_i.
\]
Let $\{p_k\}$ be orthonormal Laguerre polynomials for $w(x)=x^a e^{-x}$ and set $\varphi_k(x):=p_k(x)\,w(x)^{1/2}$. The finite-$N$ Laguerre kernel is the CD kernel
\[
K^{(L)}_N(x,y):=\sum_{k=0}^{N-1}\varphi_k(x)\varphi_k(y)
=\frac{\gamma_{N-1}}{\gamma_N}\,\frac{\varphi_N(x)\varphi_{N-1}(y)-\varphi_{N-1}(x)\varphi_N(y)}{x-y},
\]
with the usual normalization constants $\gamma_k$.

The largest-eigenvalue CDF is the one-sided gap probability
\[
F^{(L)}_N(s)=\PP(\lambda_{\max}<s)=\det\nolimits_{L^2((s,\infty))}\!\bigl(I-\Kop^{(L)}_N\bigr),
\]
where $(\Kop^{(L)}_N f)(x)=\int_s^\infty K^{(L)}_N(x,y)f(y)\,dy$.

Following \cite{forresterwitte2002pv,forrester2010log}, it is convenient to introduce the generating function
\[
E_{N,2}^{(L)}((0,t);\xi):=\det\nolimits_{L^2((0,t))}\!\bigl(I-\xi\,\Kop^{(L)}_N\bigr),
\]
for which there is a Painlevé--V $\tau$--function $\tau_V(t)$ such that
\[
E_{N,2}^{(L)}((0,t);\xi)
=t^{(a+\mu)N+N^2}\,\tau_V(t)\Big|_{\nu_0=0,\ \nu_1=-\mu,\ \nu_2=N+a,\ \nu_3=N}.
\]
Writing
\[
\sigma^{(L)}_N(t):=\frac{d}{dt}\log E^{(L)}_{N,2}\!\bigl((0,t);1\bigr),
\]
one finds that $\sigma^{(L)}_N$ solves the Jimbo--Miwa--Okamoto $\sigma$--form of $\mathrm{PV}$ with parameters $(\nu_0,\nu_1,\nu_2,\nu_3)=(0,0,N+a,N)$. The analytic behaviour at $t\to0^+$ is fixed by the LUE one-point function and singles out the random-matrix branch. At the hard edge $t\sim O(N^{-1})$ one recovers $\PIIIp$; at the soft edge, after the usual centering and $N^{1/3}$ scaling, one again finds $\PII$ and Tracy--Widom.

\subsection{Finite-\texorpdfstring{$N$}{N} JUE: Fredholm determinants and $\sigma$--\texorpdfstring{$\PVI$}{PVI}}\label{sec:PVI-theory}

For $a,b>-1$, the JUE has eigenvalue density on $(-1,1)$
\[
\frac{1}{\widetilde{Z}_{N,a,b}} \prod_{1\le i<j\le N}|\lambda_i-\lambda_j|^2
\prod_{i=1}^N (1-\lambda_i)^{a}(1+\lambda_i)^{b}\,\mathbf{1}_{(-1,1)}(\lambda_i)\,d\lambda_i.
\]
Write $\{P_k\}$ for orthonormal Jacobi polynomials for $w(x)=(1-x)^a(1+x)^b$ and set $\psi_k(x):=P_k(x)\,w(x)^{1/2}$. The finite-$N$ Jacobi kernel is
\[
K^{(J)}_N(x,y):=\sum_{k=0}^{N-1}\psi_k(x)\psi_k(y)
=\frac{\kappa_{N-1}}{\kappa_N}\,\frac{\psi_N(x)\psi_{N-1}(y)-\psi_{N-1}(x)\psi_N(y)}{x-y}.
\]
For $t\in(-1,1)$ and $\xi\in\C$,
\[
E_{N,2}^{(J)}((t,1);\xi):=\det\nolimits_{L^2((t,1))}\!\bigl(I-\xi\,\Kop^{(J)}_N\bigr),
\]
and
\[
F^{(J)}_N(s):=\PP(\lambda_{\max}<s)=E_{N,2}^{(J)}((s,1);1).
\]

There is a Painlevé--VI $\tau$--function $\tau_{VI}(t)$ whose $\sigma$--function
\[
\sigma_{VI}(t):=t(t-1)\,\frac{d}{dt}\log \tau_{VI}(t)
\]
solves the $\sigma$--form of $\PVI$ with parameters
\[
v_1=v_3=N+\frac{a+b}{2},\qquad v_2=\frac{a+b}{2},\qquad v_4=\frac{b-a}{2},
\]
see \cite{forresterwitte2004pvi,forrester2010log}. The logarithmic derivative of the JUE gap can be identified with $\sigma_{VI}$ up to an affine transform. After rescaling near the endpoints $t=\pm1$ one obtains the Bessel hard-edge limits, and in the double-Wishart (MANOVA) soft-edge regime one recovers the Tracy--Widom $F_2$ law \cite{johnstone2008jacobi,forrester2010log}.

\section{Numerical verification via Painlevé~II}

We begin with the simplest setting, where the Painlevé description is well conditioned and can serve as a reference. We solve \eqref{eq:PII} with boundary data approximating $\Ai(x)$ at large positive $x$ and then evaluate $F_2$ via \eqref{eq:TW-PII}. Independently, we compute the Fredholm determinant representation \eqref{eq:Fredholm} using a Nyström method with Gauss--Legendre quadrature.

\subsection{Numerical methods}

We briefly summarize the schemes used to evaluate the Tracy--Widom GUE distribution $F_2(x)$ on a real interval $[x_{\min},x_{\max}]$.

For the Painlevé~II side we rewrite $\PII$ as the first-order system
\[
q'(x)=p(x),\qquad p'(x)=x q(x)+2 q(x)^3
\]
for $u=(q,p)^{\mathsf T}$ and approximate the asymptotic condition at $+\infty$ by prescribing
\[
q(T_0)=\Ai(T_0),\qquad q'(T_0)=\Ai'(T_0)
\]
at some large $T_0$. We then integrate backwards from $T_0$ down to $x_{\min}$ with a high-order explicit Runge--Kutta method (in the implementation, a 7th-order Verner scheme with adaptive stepsize). The solver provides a dense representation of $q(x)$ which we sample on a uniform grid $x_k\in[x_{\min},x_{\max}]$.

Given sampled values $q(x_k)$ we approximate
\[
\log F_2(x)=-\int_x^\infty (t-x)\,q(t)^2\,dt.
\]
On the grid we compute cumulative trapezoidal sums from right to left,
\[
I_0(x_k)\approx\int_{x_k}^{x_{\max}}q(t)^2\,dt,\qquad
I_1(x_k)\approx\int_{x_k}^{x_{\max}}t\,q(t)^2\,dt,
\]
and approximate the tail beyond $T_0$ by replacing $q$ with $\Ai$ and integrating numerically. This yields
\[
\log F_2(x_k)\approx -\bigl(I_1(x_k)-x_k I_0(x_k)\bigr) - \Bigl(\int_{T_0}^{\infty} t\,\Ai(t)^2\,dt
  - x_k \int_{T_0}^{\infty} \Ai(t)^2\,dt\Bigr),
\]
and $F_2(x_k)$ is obtained by exponentiation.

On the Fredholm side we discretize the Airy kernel
\[
K_{\mathrm{Ai}}(s,t)
=\frac{\Ai(s)\Ai'(t)-\Ai'(s)\Ai(t)}{s-t},\qquad
K_{\mathrm{Ai}}(s,s)=\Ai'(s)^2-s\,\Ai(s)^2,
\]
on $(x,\infty)$ via Nyström. We map $(x,\infty)$ to $(0,1)$ by $t=x+z/(1-z)$ and pull back an $N$-point Gauss--Legendre rule $\{z_j,w_j\}$. Denoting $t_j=x+z_j/(1-z_j)$ and $\omega_j=w_j/(1-z_j)^2$, we form the symmetric matrix
\[
A_{ij}(x)=\sqrt{\omega_i}\,K_{\mathrm{Ai}}(t_i,t_j)\sqrt{\omega_j},
\]
compute the eigenvalues $\lambda_1,\dots,\lambda_N$ of $I-A(x)$, and set
\[
F_2(x)\approx\exp\Bigl(\sum_{j=1}^N\log\lambda_j\Bigr).
\]
We repeat this for the same $x$-grid as on the Painlevé side.

\subsection{Results}

\begin{figure}[ht]
    \centering
    \includegraphics[width=0.6\textwidth]{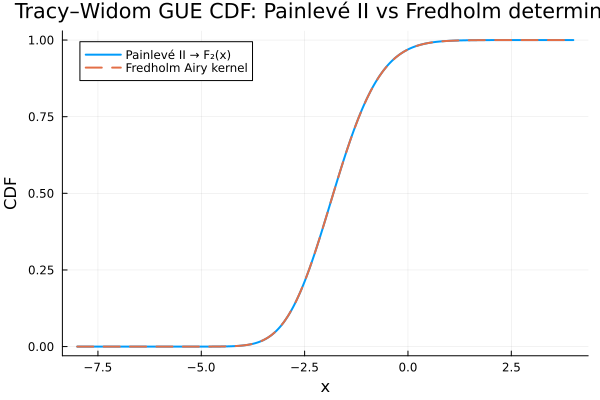}
    \caption{Tracy--Widom GUE CDF $F_2(x)$ from the Painlevé~II solution (solid) and the Airy-kernel Fredholm determinant (dashed).}
    \label{fig:cdf}
\end{figure}

\begin{figure}[ht]
    \centering
    \includegraphics[width=0.6\textwidth]{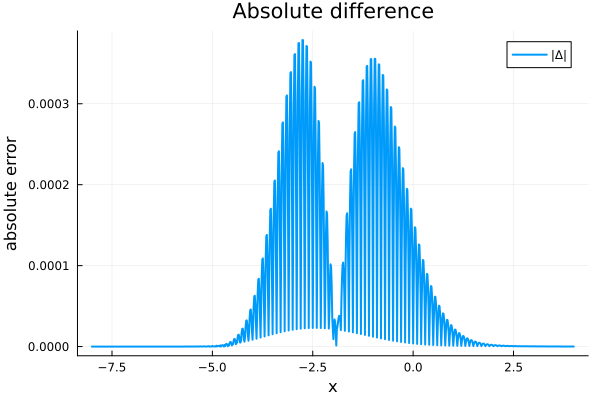}
    \caption{Absolute difference between the two numerical evaluations of $F_2(x)$.}
    \label{fig:diff}
\end{figure}

Figure~\ref{fig:cdf} shows the Tracy--Widom CDF $F_2(x)$ computed from Painlevé~II and from the Airy-kernel Fredholm determinant on $x\in[-8,4]$. The curves are visually indistinguishable. The absolute difference in Figure~\ref{fig:diff} stays below $3\times10^{-4}$, consistent with the quadrature and ODE tolerances for $N=80$ Nyström nodes.

\subsection{A brief comment on conditioning}

The Painlevé~II computation is intrinsically ill conditioned when formulated as a backward IVP from $T_0$. Linearization around the Hastings--McLeod solution reveals an unstable direction whose amplification factor between $T_0$ and the bulk is of order $\exp\bigl(c(T_0^{3/2}-x^{3/2})\bigr)$ for some $c>0$, so double-precision roundoff at $T_0$ is magnified to roughly $10^{-4}$ near $x\approx0$. Moreover we impose only the leading Airy asymptotics at $T_0$ and approximate the tails of the integral. For our purposes this four-digit accuracy is acceptable and provides a clean benchmark against which to compare the much more fragile finite-$n$ Painlevé~IV calculations. A broader discussion of numerical conditioning and limitations is given in Section~\ref{sec:limitations-future}.

\section{Numerical verification of \texorpdfstring{$\PIV$}{PIV} for finite-\texorpdfstring{$n$}{n} GUE}

We now turn to the main object of interest: the finite-$n$ GUE gap and its $\sigma$--$\PIV$ description. This is where numerical instability becomes serious and where our Fredholm-anchored $\sigma$--form integration comes into play.

\subsection{Direct IVP for the $\sigma$--form and its failure}

\begin{figure}[htbp]
  \centering
  \begin{subfigure}{0.48\textwidth}
    \imgw{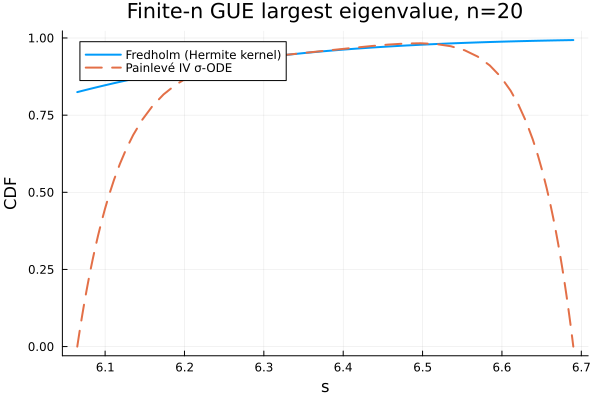}
    \caption{$n=20$: CDF comparison.}
  \end{subfigure}
  \begin{subfigure}{0.48\textwidth}
    \imgw{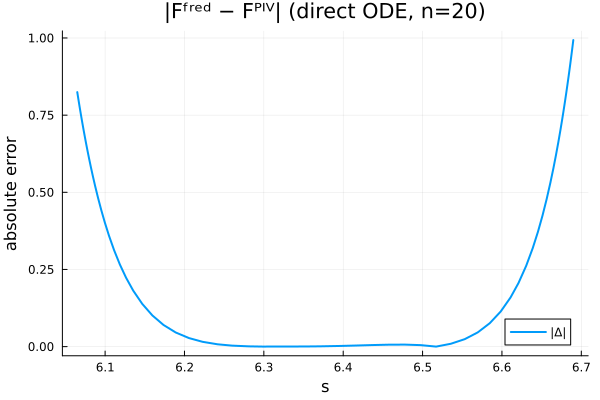}
    \caption{$n=20$: absolute error.}
  \end{subfigure}

  \vspace{0.4cm}

  \begin{subfigure}{0.48\textwidth}
    \imgw{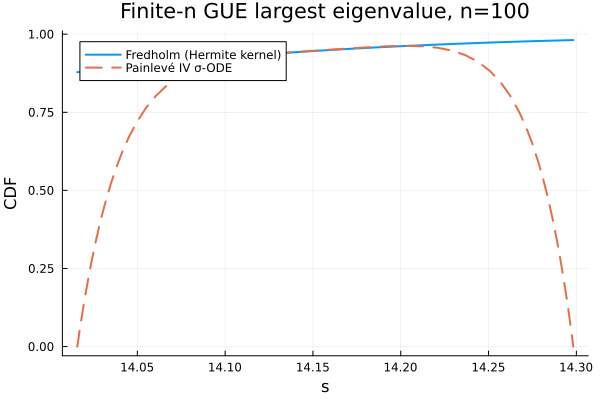}
    \caption{$n=100$: CDF comparison.}
  \end{subfigure}
  \begin{subfigure}{0.48\textwidth}
    \imgw{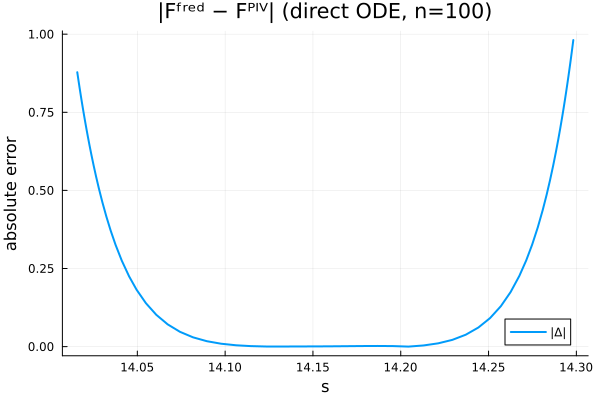}
    \caption{$n=100$: absolute error.}
  \end{subfigure}
  \caption{Direct ODE solutions of the $\PIV$ $\sigma$–equation compared with the finite-$n$ GUE Fredholm determinant.}
  \label{fig:piv-direct-ode}
\end{figure}

Figure~\ref{fig:piv-direct-ode} illustrates what happens if one simply treats \eqref{eq:PIV-sigma} as an initial-value problem. We choose a base point $s_0\approx\sqrt{2n}$, estimate $\sigma_n(s_0)$ and $\sigma_n'(s_0)$ from a local polynomial fit to $\log F_n$ computed via the Fredholm determinant, and then integrate \eqref{eq:PIV-sigma} forward and backward in $s$ with a standard ODE solver.

The $\sigma$–$\PIV$ equation is second order and second degree with movable poles, and the random-matrix solution corresponding to $F_n$ is a very special trajectory. The Cauchy data extracted at $s_0$ are only approximate, and tiny perturbations tend to push the ODE onto a different Painlevé branch that typically develops a real pole within the $s$--window of interest. Numerically this shows up as the adaptive solver repeatedly shrinking the timestep until it reaches machine precision and aborts, and as non-monotone, bump-shaped curves $F_{\mathrm{PIV}}(s)$ which match the Fredholm CDF only in a small neighbourhood of $s_0$ before bending back down toward zero at both ends.

In other words, direct shooting from a single approximate anchor is an ill-posed strategy in standard double precision: the IVP is extremely sensitive to initial data, blind to the requirement that $F_n$ be a CDF, and prone to encountering movable poles. This motivates the anchored and branch-locked approach described next.

\subsection{Numerical methods (ODE-based \texorpdfstring{$\sigma$}{sigma}-form integration)}

\emph{This subsection describes the prototype algorithm and is kept essentially unchanged, since its structure is reused in the later LUE and JUE experiments.}

We now describe a variant of the above scheme in which the $\PIV$ $\sigma$-equation is integrated between anchors using a standard ODE solver, rather than an explicit hand-crafted predictor--corrector update. The Fredholm determinant side and the extraction of anchor data are unchanged; only the propagation mechanism between anchors is modified.

\subsubsection{ODE formulation with signed branch}

Recall the $\PIV$ $\sigma$-equation for the finite-$n$ largest-eigenvalue CDF:
\begin{equation*}
  \bigl(\sigma_n''(s)\bigr)^2
  = 4\bigl(s\,\sigma_n'(s)-\sigma_n(s)\bigr)^2
    -4\bigl(\sigma_n'(s)\bigr)^2\bigl(\sigma_n'(s)+2n\bigr),
\end{equation*}
which we write concisely as
\[
  \bigl(\sigma_n''(s)\bigr)^2 = R\bigl(s,\sigma_n(s),\sigma_n'(s)\bigr),\qquad
  R(s,\sigma,\sigma')
  := 4(s\sigma'-\sigma)^2 - 4(\sigma')^2(\sigma'+2n).
\]
Since the equation is of second degree in $\sigma_n''$, it determines $\sigma_n''$ only up to a sign:
\[
  \sigma_n''(s) = \pm \sqrt{R\bigl(s,\sigma_n(s),\sigma_n'(s)\bigr)}.
\]

To integrate numerically we rewrite the problem as a first-order system in the state vector $u=(\sigma_n,\sigma_n')$:
\[
  \frac{d}{ds}
  \begin{pmatrix}
    \sigma_n(s) \\[0.2em]
    \sigma_n'(s)
  \end{pmatrix}
  =
  \begin{pmatrix}
    \sigma_n'(s) \\[0.2em]
    \sigma_n''(s)
  \end{pmatrix}
  =
  \begin{pmatrix}
    \sigma_n'(s) \\[0.2em]
    \mathrm{sgn}_j \sqrt{R\bigl(s,\sigma_n(s),\sigma_n'(s)\bigr)}
  \end{pmatrix},
\]
where $\mathrm{sgn}_j\in\{-1,+1\}$ is a fixed branch-sign parameter on each anchor interval $[s^{(j+1)},s^{(j)}]$ (see below). Thus, once $\mathrm{sgn}_j$ is chosen, the right-hand side is single-valued and we obtain a well-posed first-order ODE system $u'(s) = f\bigl(s,u(s);n,\mathrm{sgn}_j\bigr)$, which we solve with a high-order explicit Runge--Kutta method with adaptive stepsize control.

For numerical robustness we include a safeguard when the radicand becomes non-positive: if $R(s,\sigma,\sigma')\le 0$ we set the second component of $f(s,u)$ to zero, effectively freezing $\sigma_n''$ at that point. In practice this only occurs deep in the tails, where $F_n(s)$ is extremely close to $0$ or $1$ and the overall error is dominated by the Fredholm truncation.

\subsubsection{Anchor data and branch selection on each interval}

We fix a spectral window
\[
  [s_{\min},s_{\max}]
  = [\sqrt{2n}-W,\sqrt{2n}+W],
\]
with $W\approx 3$, and choose a descending sequence of anchor points $s^{(1)} > s^{(2)} > \cdots > s^{(N_{\mathrm{anc}})}$ uniformly in $[s_{\min},s_{\max}]$. At each anchor $s^{(j)}$ we compute $F_n(s)$ on a short symmetric stencil and fit a degree-4 polynomial to
\[
  g(s) := \log F_n(s).
\]
Writing
\[
  g(s^{(j)}+u) \approx
    c_0^{(j)} + c_1^{(j)}u + c_2^{(j)}u^2
    + c_3^{(j)}u^3 + c_4^{(j)}u^4,
\]
we recover the approximate Cauchy data
\[
  F_n(s^{(j)}) \approx e^{c_0^{(j)}},\qquad
  \sigma_n(s^{(j)}) \approx c_1^{(j)},\qquad
  \sigma_n'(s^{(j)}) \approx 2c_2^{(j)},\qquad
  \sigma_n''(s^{(j)}) \approx 6c_3^{(j)}.
\]
These values are then used to define both the initial condition and the branch sign on the interval to the left of $s^{(j)}$.

For each interval $[s^{(j+1)},s^{(j)}]$ we set
\[
  u^{(j)}(s^{(j)})
    := \bigl(\sigma_n(s^{(j)}),\sigma_n'(s^{(j)})\bigr)
     \approx \bigl(c_1^{(j)},2c_2^{(j)}\bigr),
\]
and choose
\[
  \mathrm{sgn}_j :=
  \begin{cases}
    +1, & \text{if } \bigl|\sigma_n''(s^{(j)})\bigr|
                      \le \varepsilon_{\mathrm{init}},\\[0.3em]
    \mathrm{sign}\bigl(\sigma_n''(s^{(j)})\bigr),
      & \text{otherwise},
  \end{cases}
\]
where $\sigma_n''(s^{(j)})\approx 6c_3^{(j)}$ and $\varepsilon_{\mathrm{init}}>0$ is a small threshold. This ensures that each interval carries a constant branch sign tied to the Fredholm-based estimate of $\sigma_n''$ at its right endpoint while preventing the integration from being trapped near $\sigma_n''\equiv 0$ when the third derivative is numerically small.

\subsubsection{ODE integration between anchors and reprojection}

We introduce a global descending grid
\[
  s_1 = s_{\max} > s_2 > \cdots > s_{N_{\mathrm{grid}}} = s_{\min}
\]
covering the same spectral window and partition it by anchor intervals. On each interval $[s^{(j+1)},s^{(j)}]$ we form the list of grid points with $s_k\in[s^{(j+1)},s^{(j)}]$ and integrate
\[
  u'(s) = f\bigl(s,u(s);n,\mathrm{sgn}_j\bigr),
\]
with initial condition $u(s^{(j)}) = u^{(j)}(s^{(j)})$, evaluating the solution at these grid points. This yields values $(\sigma_n(s_k),\sigma_n'(s_k))$ with local error controlled by the usual ODE tolerances.

At the left endpoint $s^{(j+1)}$ we re-anchor the solution by overwriting the propagated values with the Fredholm-based Cauchy data:
\[
  \sigma_n(s^{(j+1)}) \leftarrow c_1^{(j+1)},\qquad
  \sigma_n'(s^{(j+1)}) \leftarrow 2c_2^{(j+1)}.
\]
The corresponding $\sigma_n''(s^{(j+1)})$ is retained only for determining the next branch sign. This reprojection prevents the accumulation of drift and keeps the ODE solution consistent with the Hermite-kernel Fredholm determinant at every anchor.

\subsubsection{Reconstruction of the CDF}

Once $\sigma_n(s_k)$ has been computed on the global grid, we reconstruct $F_n(s_k)$ from
\[
\frac{d}{ds}\log F_n(s)=\sigma_n(s)
\]
by trapezoidal integration from a base anchor $s_0$ with known $F_n(s_0)$:
\[
  \log F_n(s_k) \approx \log F_n(s_0)
   + \sum_{\ell:\,s_\ell\in[s_0,s_k]}
     \frac{\sigma_n(s_\ell)+\sigma_n(s_{\ell+1})}{2}\,(s_{\ell+1}-s_\ell).
\]
We then set $F_n(s_k)=\exp(\log F_n(s_k))$. This gives a numerical CDF which agrees with the Hermite Fredholm determinant to high accuracy and in which the $\PIV$ $\sigma$--form acts as a dynamical constraint between sparse Fredholm anchors.

\subsection{Results and analysis}

\begin{figure}[htbp]
\centering
\begin{subfigure}{0.48\textwidth}
\imgw{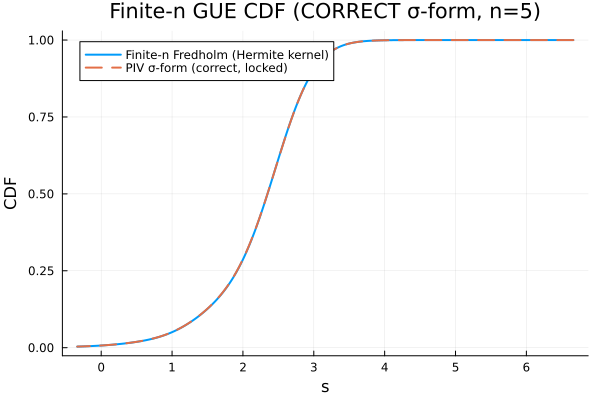}
\caption{$n=5$: CDF.}
\label{fig:cdf-n5}
\end{subfigure}
\begin{subfigure}{0.48\textwidth}
\imgw{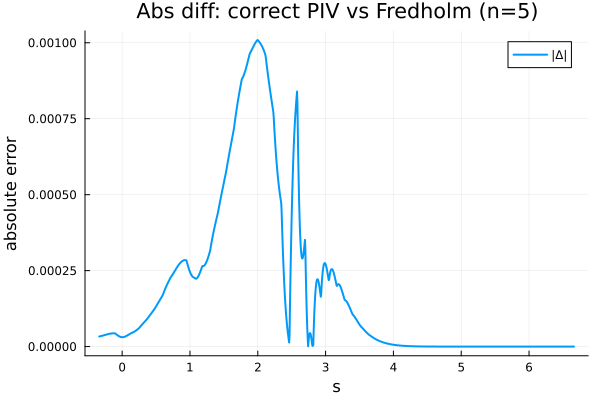}
\caption{$n=5$: absolute error.}
\label{fig:error-n5}
\end{subfigure}
\end{figure}

\begin{figure}[htbp]
\centering
\begin{subfigure}{0.48\textwidth}
\imgw{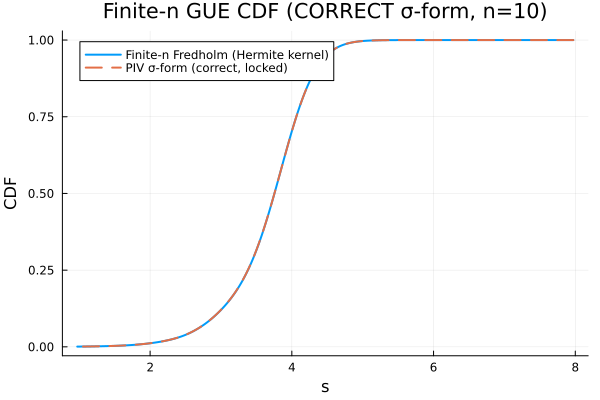}
\caption{$n=10$: CDF.}
\label{fig:cdf-n10}
\end{subfigure}
\begin{subfigure}{0.48\textwidth}
\imgw{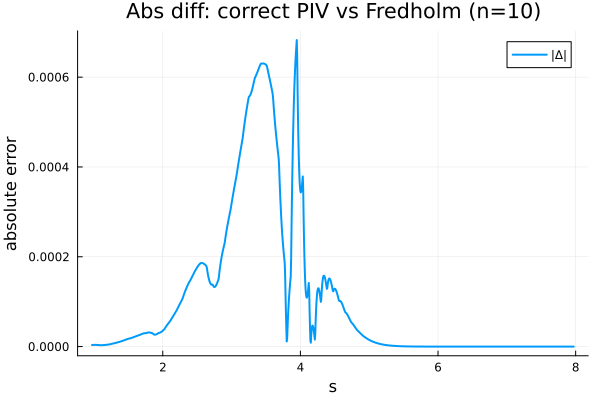}
\caption{$n=10$: absolute error.}
\label{fig:error-n10}
\end{subfigure}
\end{figure}

\begin{figure}[htbp]
\centering
\begin{subfigure}{0.48\textwidth}
\imgw{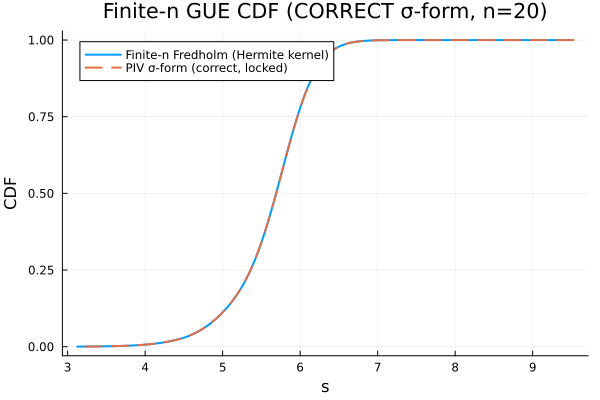}
\caption{$n=20$: CDF.}
\label{fig:cdf-n20}
\end{subfigure}
\begin{subfigure}{0.48\textwidth}
\imgw{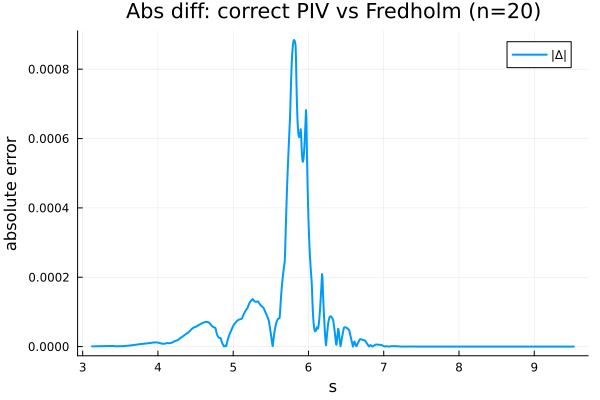}
\caption{$n=20$: absolute error.}
\label{fig:error-n20}
\end{subfigure}
\end{figure}

\begin{figure}[htbp]
\centering
\begin{subfigure}{0.48\textwidth}
\imgw{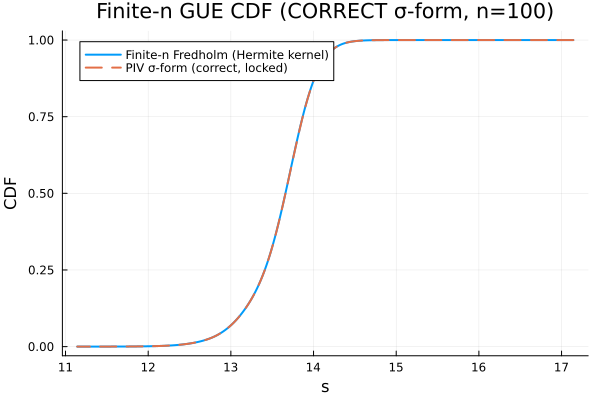}
\caption{$n=100$: CDF.}
\label{fig:cdf-n100}
\end{subfigure}
\begin{subfigure}{0.48\textwidth}
\imgw{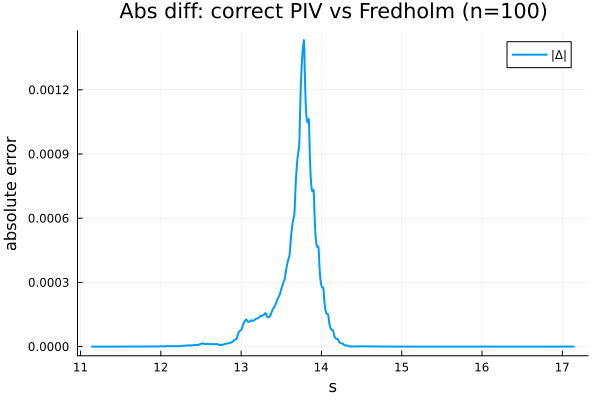}
\caption{$n=100$: absolute error.}
\label{fig:error-n100}
\end{subfigure}
\end{figure}

\begin{figure}[htbp]
\centering
\begin{subfigure}{0.48\textwidth}
\imgw{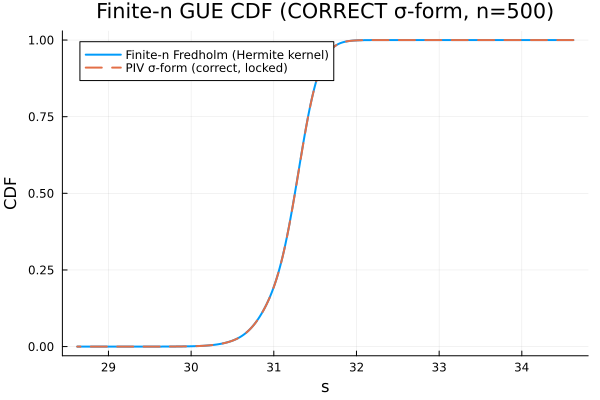}
\caption{$n=500$: CDF.}
\label{fig:cdf-n500}
\end{subfigure}
\begin{subfigure}{0.48\textwidth}
\imgw{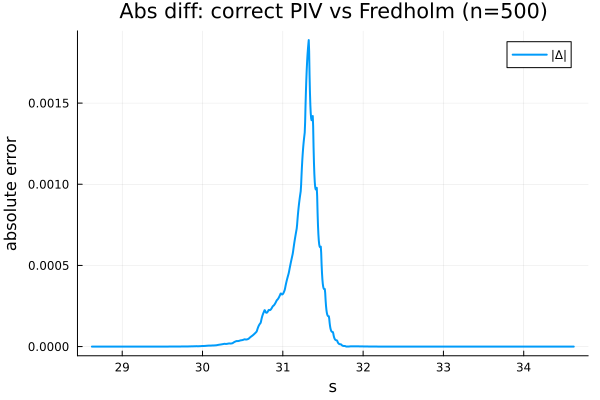}
\caption{$n=500$: absolute error.}
\label{fig:error-n500}
\end{subfigure}
\end{figure}

Figures~\ref{fig:cdf-n5}--\ref{fig:error-n500} show the comparison between the anchored $\PIV$ solution and the Hermite Fredholm determinant for $n=5,10,20,100,500$. In each case we display the CDF and the absolute error over a window centered at $\sqrt{2n}$.

Table~\ref{tab:errors} summarizes the maximum absolute errors between the two methods and the corresponding spectral windows.

\begin{table}[h]
\centering
\begin{tabular}{@{}cccc@{}}
\toprule
$n$ & Max error & Edge location & Window width \\
\midrule
  5   & $1.01 \times 10^{-3}$ & $\sqrt{10}\approx3.16$   & 7.0 \\
 10   & $6.83 \times 10^{-4}$ & $\sqrt{20}\approx4.47$   & 7.0 \\
 20   & $8.84 \times 10^{-4}$ & $\sqrt{40}\approx6.32$   & 6.4 \\
100   & $1.43 \times 10^{-3}$ & $\sqrt{200}\approx14.14$ & 6.0 \\
500   & $1.89 \times 10^{-3}$ & $\sqrt{1000}\approx31.62$& 6.0 \\
\bottomrule
\end{tabular}
\caption{Maximum absolute errors and problem parameters for the anchored $\PIV$ method.}
\label{tab:errors}
\end{table}

For $n$ up to about $500$ the anchored $\PIV$ and Fredholm CDFs agree to roughly $10^{-3}$ in the transition region and much better in the tails. For very large $n$ the fixed physical window $[s_{\min},s_{\max}]$ sits far into the regime where $F_n(s)$ is already extremely close to $1$, and both the Hermite basis and the Fredholm matrix become ill conditioned in double precision. At that point the true error is dominated by underflow and cancellation, and the plots primarily reflect saturation effects rather than meaningful extra digits. Conceptually this is not a problem: for such $n$ it is more natural to pass to the soft-edge scaling and use the stable $\PII$/$F_2$ description. A more systematic discussion of the global limitations of the anchored scheme is deferred to Section~\ref{sec:limitations-future}.

\subsection{Large-$n$ behaviour and interpretation}

For moderate sizes the anchored $\PIV$ computation behaves exactly as dictated by the theory: it recovers the finite-$n$ GUE gap with a controlled discrepancy that is largest near the soft edge and negligible in the tails. As $n$ increases, the natural soft-edge scale is
\[
s=\mu(n)+\sigma(n)x,\qquad \mu(n)=\sqrt{2n},\quad \sigma(n)\asymp n^{-1/6},
\]
so a fixed physical window $[\,\mu(n)-\Delta,\mu(n)+\Delta\,]$ corresponds to an $x$-window that grows with $n$. On the right half of this window the Tracy--Widom limit $F_2(x)$ is already extremely close to $1$ when $x$ is moderately positive, and double-precision arithmetic can no longer resolve small differences in $\log F_n$ or in the Nyström eigenvalues. The anchored $\PIV$ solver then essentially returns a plateau at $F_n\approx 1$ with tiny numerical fluctuations. For the purposes of this paper we do not push beyond $n\approx500$ in the finite-$n$ GUE computations and treat the large-$n$ regime via $\PII$ and the Airy kernel instead.

\section{Painlevé~VI for JUE and finite-\texorpdfstring{$N$}{N} hard edges}

\subsection{Finite-\texorpdfstring{$N$}{N} JUE gap and anchored \texorpdfstring{$\PVI$}{PVI}}\label{sec:num-JUE-hard}

For the JUE we verify that the largest-eigenvalue CDF
\[
F_N(s)=\det\nolimits_{L^2(s,1)}(I-K_N)
\]
agrees, for finite $N$, with the Painlevé~VI description derived from the Jimbo--Miwa--Okamoto $\sigma$--form. The finite-$N$ kernel $K_N$ is constructed from orthonormal Jacobi functions $\phi_k$ as in Section~\ref{sec:PVI-theory}. To approximate $F_N(s)$ we use a Gram-type Nyström method: map Gauss--Legendre nodes $z_i\in(-1,1)$ to $t_i\in(s,1)$, form the matrix $\Phi_{ik}=\sqrt{w_i^{(s)}}\,\phi_k(t_i)$, and compute
\[
F_N^{\mathrm{(Fred)}}(s)\approx \det(I_N-G),\qquad G=\Phi^\mathsf T\Phi.
\]

On the Painlevé side we adopt the same philosophy as for $\PIV$. We work with the $\PVI$ $\sigma$--form associated to the JUE gap and the corresponding second-order second-degree ODE for $f(t)=F_N^{(J)}(t;a,b,0;1)$. Introducing $g(t)=(1-t)f'(t)$ and the stretched variable $y=-\log(1-t)$, we obtain a first-order system for $(f,g,t)$ in $y$, with $f''(t)$ recovered at each step from the algebraic $\sigma$--form. The two quadratic roots are disambiguated by a branch-locking rule exactly as in the $\PIV$ case.

Anchors are now placed automatically in the $y$-variable: for each $(N,a,b)$ we determine a window $[s_{\min},s_{\max}]$ by solving $F_N^{\mathrm{(Fred)}}(s_{\min})\approx 0.1$ and $F_N^{\mathrm{(Fred)}}(s_{\max})\approx 1-10^{-8}$, then map this to $t\in[t_{\min},t_{\max}]$ and $y\in[y_{\min},y_{\max}]$. Uniformly spaced $y_j$ are mapped back to $s_j$ and used as anchors. At each $s_j$ we extract $F_N(s_j)$ and its first two derivatives in $s$ from a local least-squares fit to $\log F_N(s)$, convert these to $f(t_j)$, $f'(t_j)$, $f''(t_j)$, and launch a PVI integration segment in $y$ until the next anchor, where we relock the solution.

As in the GUE case we find that the anchored $\PVI$ solution reproduces the JUE Fredholm CDF with errors of size $10^{-3}$--$10^{-2}$ in the transition region, relatively insensitive to $N$ once the quadrature and anchor grids are chosen. Rather than repeating a detailed error discussion here, we refer to the hard-edge analysis below and to the global limitations in Section~\ref{sec:limitations-future}.

\section{JUE hard-edge limits and Bessel kernels}\label{sec:num-JUE-hard-limits}

We next check numerically that the Jacobi ensemble with weight $(1-x)^a(1+x)^b$ has Bessel hard-edge limits at both endpoints. Writing $F_N^{(J)}(t;a,b)$ for the finite-$N$ gap probability, the theory predicts
\[
E^{(a,b)}_{N,\mathrm{right}}(s)
  := F_N^{(J)}\Bigl(1-\frac{s}{2N^2};a,b\Bigr)\longrightarrow E_{\mathrm{hard}}^{(a)}(s),
\]
\[
E^{(a,b)}_{N,\mathrm{left}}(s)
  := F_N^{(J)}\Bigl(-1+\frac{s}{2N^2};a,b\Bigr)
   = F_N^{(J)}\Bigl(1-\frac{s}{2N^2};b,a\Bigr)\longrightarrow E_{\mathrm{hard}}^{(b)}(s),
\]
where $E_{\mathrm{hard}}^{(\alpha)}(s)=\det(I-K^{(\alpha)}_{\mathrm{Bes}})_{L^2(0,s)}$ and $K^{(\alpha)}_{\mathrm{Bes}}$ is the Bessel kernel of order $\alpha$ \cite{tracy1994bessel,deift1999orthogonal,forrester2010log}.

On the finite-$N$ side we evaluate $F_N^{(J)}$ by Nyström discretization of the rank-$N$ Jacobi kernel as before, now on $(t,1)$ with $t$ near $\pm1$ and $t$ scaled as above. On the Bessel side we implement the standard off-diagonal formula
\[
K^{(\alpha)}_{\mathrm{Bes}}(x,y)
 =
 \frac{J_\alpha(\sqrt{x}) \sqrt{y}\,J'_\alpha(\sqrt{y})
     - J_\alpha(\sqrt{y}) \sqrt{x}\,J'_\alpha(\sqrt{x})}
      {2(x-y)},
\]
with the diagonal handled by a symmetric finite difference, and discretize on $(0,s)$ using Gauss--Legendre nodes. In both cases the eigenvalues of $I-A$ are used to compute $\log\det(I-A)$.

For the symmetric case $(a,b)=(0,0)$ the left and right edges coincide. Table~\ref{tab:hardedge-00} reports the maximal absolute discrepancy over $s\in[0,15]$ for $N=20,40,80,120,300$.

\begin{table}[ht]
  \centering
  \caption{Maximum hard-edge error for $(a,b)=(0,0)$:
  $\max_s |E^{(0,0)}_{N,\mathrm{edge}}(s)-E^{(0)}_{\mathrm{hard}}(s)|$ on $s\in[0,15]$.}
  \label{tab:hardedge-00}
  \begin{tabular}{@{}ccc@{}}
    \toprule
    $N$ & right edge max $|\Delta|$ & left edge max $|\Delta|$ \\
    \midrule
     20  & $6.77\times 10^{-4}$ & $6.77\times 10^{-4}$ \\
     40  & $1.69\times 10^{-4}$ & $1.69\times 10^{-4}$ \\
     80  & $4.23\times 10^{-5}$ & $4.23\times 10^{-5}$ \\
    120  & $1.88\times 10^{-5}$ & $1.88\times 10^{-5}$ \\
    300  & $4.25\times 10^{-6}$ & $4.25\times 10^{-6}$ \\
    \bottomrule
  \end{tabular}
\end{table}

The error decreases rapidly with $N$ and is already at the $10^{-6}$ level for $N=300$. For asymmetric parameters the convergence is slightly slower but still robust. Table~\ref{tab:hardedge-300} lists the maximal errors for $N=300$ and several $(a,b)$.

\begin{table}[ht]
  \centering
  \caption{Maximum hard-edge error at $N=300$:
  $\max_s |E^{(a,b)}_{N,\mathrm{edge}}(s)-E^{(\alpha)}_{\mathrm{hard}}(s)|$ on $s\in[0,15]$, where $\alpha=a$ at the right edge and $\alpha=b$ at the left edge.}
  \label{tab:hardedge-300}
  \begin{tabular}{@{}cccc@{}}
    \toprule
    $(a,b)$ & edge & max $|\Delta|$ & Bessel order $\alpha$ \\ \midrule
    $(0,0)$ & right & $4.25\times 10^{-6}$ & $0$ \\
            & left  & $4.25\times 10^{-6}$ & $0$ \\[0.2em]
    $(2,0)$ & right & $4.03\times 10^{-3}$ & $2$ \\
            & left  & $2.45\times 10^{-3}$ & $0$ \\[0.2em]
    $(0,3)$ & right & $3.66\times 10^{-3}$ & $0$ \\
            & left  & $2.88\times 10^{-3}$ & $3$ \\[0.2em]
    $(2,3)$ & right & $1.01\times 10^{-2}$ & $2$ \\
            & left  & $4.82\times 10^{-3}$ & $3$ \\
    \bottomrule
  \end{tabular}
\end{table}

\begin{figure}[ht]
  \centering
  \begin{minipage}{0.48\textwidth}
    \includegraphics[width=\linewidth]{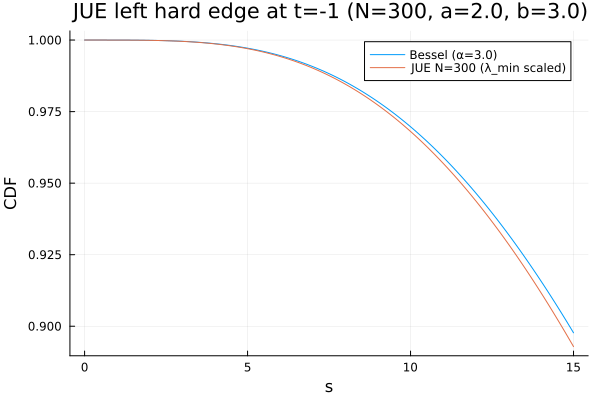}
    
    \includegraphics[width=\linewidth]{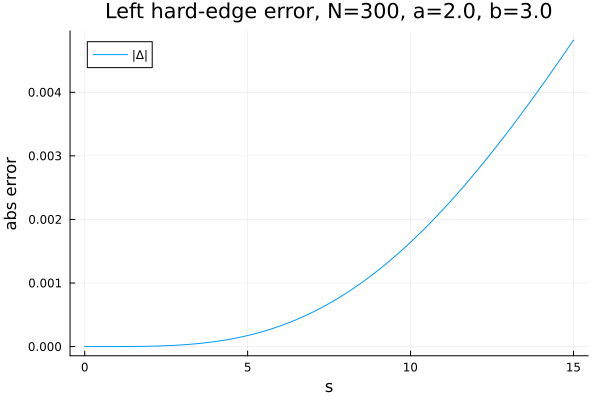}
  \end{minipage}\hfill
  \begin{minipage}{0.48\textwidth}
    \includegraphics[width=\linewidth]{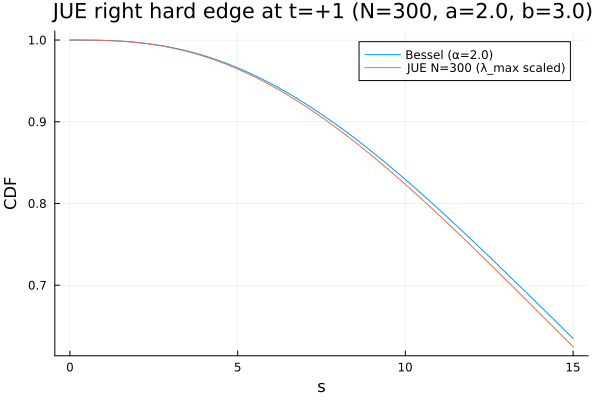}
    
    \includegraphics[width=\linewidth]{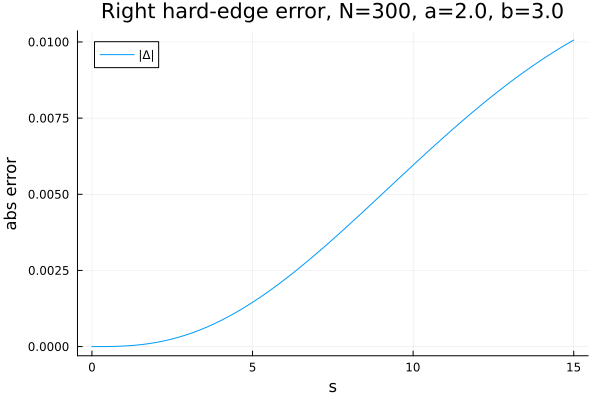}
  \end{minipage}
  \caption{Hard-edge comparison for $(a,b)=(2,3)$ at $N=300$. Top: CDF of rescaled JUE extremal eigenvalues versus the Bessel hard-edge limit at the left (top left) and right (top right) edges. Bottom: corresponding absolute errors.}
  \label{fig:jue-hardedge-23}
\end{figure}

Figure~\ref{fig:jue-hardedge-23} shows representative CDFs and error curves for the most asymmetric case $(a,b)=(2,3)$ at $N=300$. The finite-$N$ JUE hard-edge statistics closely track the corresponding Bessel limits, with absolute errors of order $10^{-3}$--$10^{-2}$ concentrated near the edge.

\section{JUE soft edge and Tracy--Widom limit}

Finally we test the JUE soft-edge scaling limit by Monte Carlo simulation of double-Wishart (MANOVA) matrices and comparison with the standardised Tracy--Widom $F_2$ law.

We simulate the largest eigenvalue of 
\[
\Theta=(A+B)^{-1}B,
\]
where $A\sim\mathrm{CW}_N(I,n_1)$ and $B\sim\mathrm{CW}_N(I,n_2)$ are independent complex Wishart matrices with $n_1\approx2N$, $n_2\approx3N$. The eigenvalues of $\Theta$ form a JUE on $(0,1)$ with parameters determined by $(N,n_1,n_2)$, and $\lambda_{\max}$ represents the JUE soft edge. For each $N\in\{100,200,400\}$ we generate $M$ independent samples, estimate the empirical mean $\mu_N$ and standard deviation $\sigma_N$ of $\lambda_{\max}$, and form standardised variables
\[
z_m=\frac{\lambda_{\max}^{(m)}-\mu_N}{\sigma_N}.
\]
The empirical CDF $F_N(z)$ is then compared on $[-4,4]$ with the standardised Tracy--Widom CDF
\[
F_2^{\mathrm{std}}(z):=F_2\bigl(\mu_{F_2}+\sigma_{F_2}z\bigr),
\]
where $\mu_{F_2}$ and $\sigma_{F_2}$ are the numerically computed mean and standard deviation of $F_2$ from the Fredholm Airy kernel.

Table~\ref{tab:JUE_soft_edge_MC} summarises the results, including the scale-invariant quantity $\sigma_N N^{2/3}$ and the Kolmogorov distance
\[
\Delta_N:=\sup_{z\in[-4,4]}|F_N(z)-F_2^{\mathrm{std}}(z)|.
\]

\begin{table}[t]
  \centering
  \begin{tabular}{@{}cccccc@{}}
    \toprule
    $N$ & $M$ & $\mathbb E[\lambda_{\max}]$ & $\mathrm{sd}(\lambda_{\max})$
        & $\sigma_N N^{2/3}$ & $\Delta_N$ on $[-4,4]$ \\
    \midrule
    100 & 2000 & $0.94534$ & $0.00346$ & $0.07451$ & $0.017$ \\
    200 & 2000 & $0.94793$ & $0.00211$ & $0.07231$ & $0.014$ \\
    400 & 2000 & $0.94941$ & $0.00125$ & $0.06798$ & $0.016$ \\
    \midrule
    100 & 5000 & $0.94537$ & $0.00348$ & $0.07489$ & $0.018$ \\
    200 & 5000 & $0.94782$ & $0.00214$ & $0.07316$ & $0.008$ \\
    400 & 5000 & $0.94939$ & $0.00128$ & $0.06943$ & $0.012$ \\
    \bottomrule
  \end{tabular}
  \caption{Monte Carlo approximation of the JUE soft-edge law for $(N,n_1,n_2)\approx(N,2N,3N)$, compared with the standardised Tracy--Widom $F_2^{\mathrm{std}}(z)$ on $[-4,4]$.}
  \label{tab:JUE_soft_edge_MC}
\end{table}

The products $\sigma_N N^{2/3}$ are slowly varying and remain of order $10^{-1}$ as $N$ grows, consistent with the $N^{-2/3}$ soft-edge scale. The Kolmogorov distances $\Delta_N$ lie in the range $10^{-2}$--$2\times10^{-2}$, which is compatible with the expected $O(N^{-2/3})$ finite-$N$ bias and the $O(M^{-1/2})$ Monte Carlo fluctuations for $M\le 5000$. In other words, at the level of accuracy accessible with these sample sizes the JUE soft edge behaves as predicted by the Tracy--Widom $F_2$ law.

\section{LUE soft and hard edges via \texorpdfstring{$\mathrm{PV}$}{PV}}

\subsection{Finite-\texorpdfstring{$N$}{N} LUE and anchored $\sigma$--\texorpdfstring{$\mathrm{PV}$}{PV}}

For LUE we repeat the anchored $\sigma$--form strategy but now for the $\sigma$--$\mathrm{PV}$ equation described in Section~\ref{sec:LUE-PV}. The orthonormal Laguerre functions are
\[
\phi_k(x) = \frac{L_k^{(\alpha)}(x)}{\sqrt{h_k}}\,x^{\alpha/2}e^{-x/2},\qquad
h_k=\frac{\Gamma(k+\alpha+1)}{k!},
\]
and the CD kernel $K_N$ is built from $\phi_k$ with a diagonal fallback on $x=y$. We discretise $(s,\infty)$ via Nyström on $[s,s+L]$ with $M$ Gauss--Legendre nodes and form the symmetric matrix $A$ as before. The CDF $F_N(s)$ is again computed as $\exp\sum_j\log\lambda_j$ from the eigenvalues of $I-A$.

Anchors are placed near the soft edge, and local polynomial fits to $\log F_N$ provide estimates of $\sigma$ and its first derivative; these fix the branch of the algebraic $\sigma$--$\mathrm{PV}$ equation for $\sigma''$ on each interval. The ODE between anchors is integrated with the same signed-branch and reprojection strategy as in the $\PIV$ case.

\begin{figure}[t]
  \centering
  \begin{subfigure}{0.48\textwidth}
    \imgw{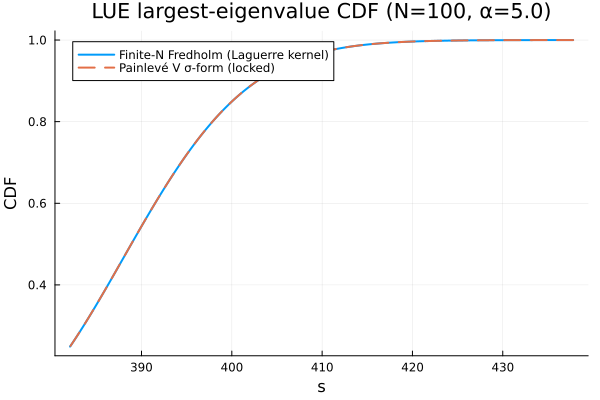}
    \caption{$N=100,\ \alpha=5$ (CDF).}
  \end{subfigure}\hfill
  \begin{subfigure}{0.48\textwidth}
    \imgw{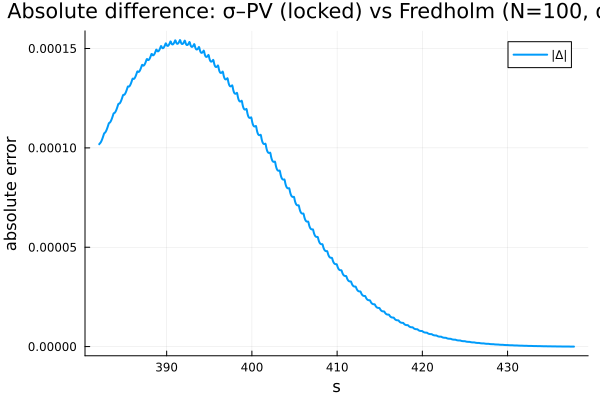}
    \caption{Absolute error.}
  \end{subfigure}
  \caption{LUE: Fredholm vs.\ locked $\sigma$--$\mathrm{PV}$ (case 1).}
  \label{fig:lue-100}
\end{figure}

\begin{figure}[t]
  \centering
  \begin{subfigure}{0.48\textwidth}
    \imgw{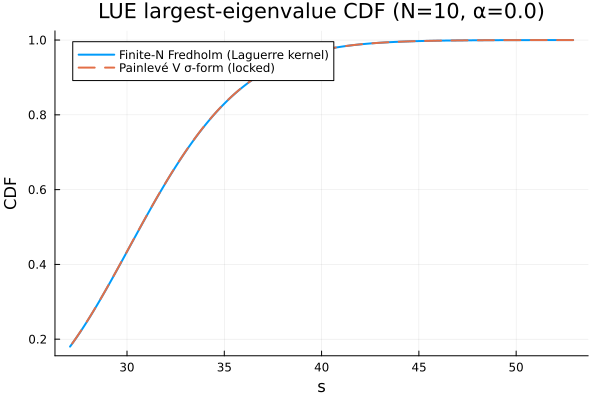}
    \caption{$N=10,\ \alpha=0$ (CDF).}
  \end{subfigure}\hfill
  \begin{subfigure}{0.48\textwidth}
    \imgw{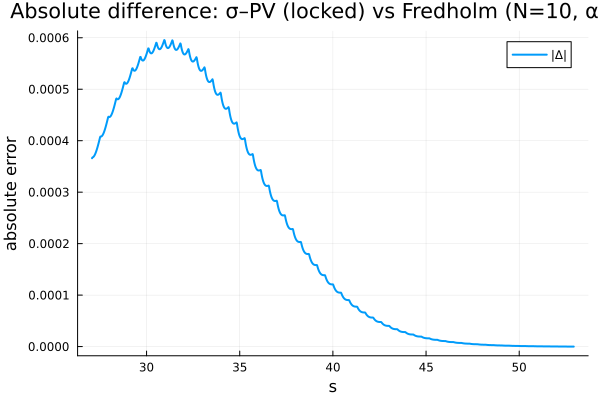}
    \caption{Absolute error.}
  \end{subfigure}
  \caption{LUE: Fredholm vs.\ locked $\sigma$--$\mathrm{PV}$ (case 4).}
  \label{fig:lue-10}
\end{figure}

Figures~\ref{fig:lue-100}--\ref{fig:lue-10} show representative comparisons between the LUE Fredholm CDF and the anchored $\sigma$--$\mathrm{PV}$ solution for four parameter choices. Table~\ref{tab:maxerr-lue} summarises the maximum absolute discrepancy, which remains below $6\times10^{-4}$ even for the smallest $N$ considered.

\begin{figure}[t]
  \centering
  \begin{subfigure}{0.48\textwidth}
    \imgw{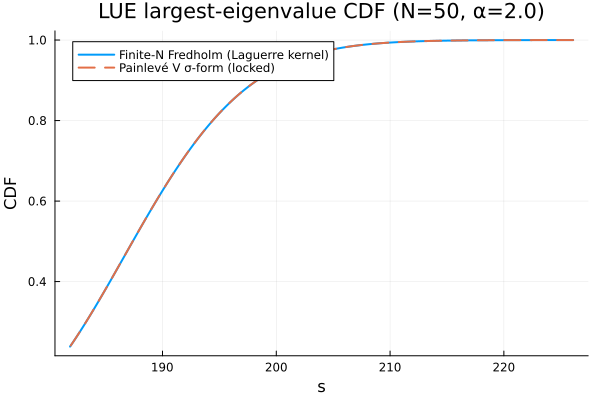}
    \caption{$N=50,\ \alpha=2$ (CDF).}
  \end{subfigure}\hfill
  \begin{subfigure}{0.48\textwidth}
    \imgw{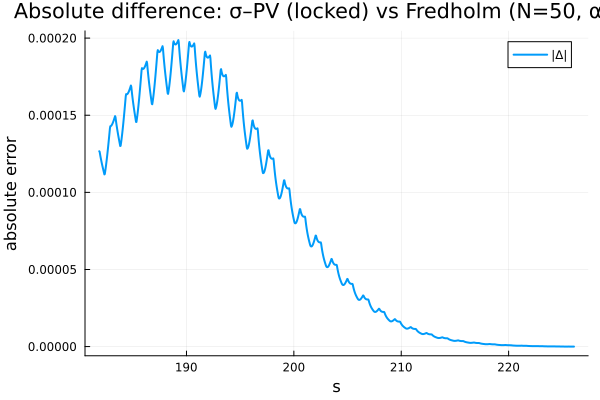}
    \caption{Absolute error.}
  \end{subfigure}
\end{figure}

\begin{figure}[t]
  \centering
  \begin{subfigure}{0.48\textwidth}
    \imgw{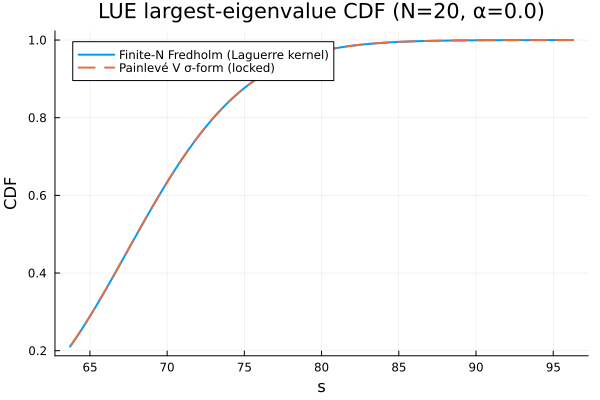}
    \caption{$N=20,\ \alpha=0$ (CDF).}
  \end{subfigure}\hfill
  \begin{subfigure}{0.48\textwidth}
    \imgw{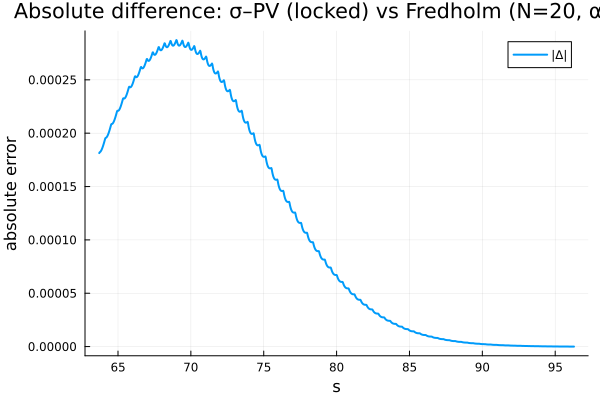}
    \caption{Absolute error.}
  \end{subfigure}
\end{figure}

\begin{table}[h]
  \centering
  \caption{Maximum absolute error $\max_s |F_N^{(\mathrm{Fred})}(s)-F_N^{(\sigma\text{-PV})}(s)|$
           and run parameters (LUE).}
  \label{tab:maxerr-lue}
  \vspace{0.5ex}
  \begin{tabular}{@{}rcccccc@{}}
    \toprule
    $N$ & $\alpha$ & $\max|\Delta|$ & $M$ (GL) & \#anchors & $s_{\min}$ & $s_{\max}$ \\
    \midrule
    100 & 5  & $1.54\times10^{-4}$ & 280 & 101 & 382.089 & 437.789 \\
     50 & 2  & $1.99\times10^{-4}$ & 260 &  91 & 181.876 & 226.085 \\
     20 & 0  & $2.87\times10^{-4}$ & 240 &  81 &  63.713 &  96.287 \\
     10 & 0  & $5.95\times10^{-4}$ & 220 &  61 &  27.073 &  52.927 \\
    \bottomrule
  \end{tabular}
\end{table}

\subsection{LUE hard and soft edges}

The same finite-$N$ Laguerre Fredholm machinery also allows us to reproduce the universal LUE hard-edge and soft-edge limits. At the hard edge we consider
\[
E_N^{(L)}((0,t);1)
  = \det\nolimits_{L^2((0,t))}\!\bigl(I-K_N^{(L)}\bigr)
\]
with $a=0$ and scaling $t=s/(4N)$, and compare it to the Bessel gap $E_{\mathrm{hard}}^{(0)}(s)$. The implementation is identical in spirit to the JUE hard-edge test: we discretise the Laguerre and Bessel kernels on $(0,t)$ and $(0,s)$, respectively, with Gauss--Legendre nodes and compare the determinants. For $N=20,40,80$ we find maximal errors of $O(10^{-4})$, $O(10^{-5})$, and $O(10^{-5})$ respectively over $s\in[0.5,10]$.

At the soft edge we study the largest-eigenvalue CDF $F_N^{(L)}(s)$ and compare its centred and scaled version with the Tracy--Widom $F_2$ law, now using fully deterministic Fredholm computations on both sides. We calibrate the centering and scaling parameters $(\mu_N,\sigma_N)$ for $N=500$ and $\alpha=0$ by matching three quantiles between $F_N^{(L)}$ and $F_2$. This yields $\mu_{500}\approx1999.86$ and $\sigma_{500}\approx19.85$, very close to the expected edge and $N^{1/3}$ scale. On $x\in[-6,6]$ the scaled LUE CDF $F_{500}^{(L)}(\mu_{500}+\sigma_{500}x)$ and $F_2(x)$ then differ by at most about $2.6\times10^{-3}$, with typical errors near $10^{-3}$ in the central region.

\section{Global limitations, Hamiltonian attempts, and future directions}\label{sec:limitations-future}

The experiments above show that, when heavily constrained by Fredholm data, Painlevé $\sigma$--forms can be used to reconstruct finite-$n$ gap probabilities with $10^{-3}$--$10^{-5}$ accuracy across GUE, LUE, and JUE. At the same time they reveal several structural limitations of this approach and point towards more robust analytic tools.

\subsection{Anchors, near interpolation, and sensitivity to perturbations}

The core strategy of the paper is to view the Painlevé $\sigma$--equations as dynamical constraints between a relatively dense set of Fredholm anchors. At each anchor we estimate $F$, $\sigma=\dv{}{s}\log F$, and $\sigma'$ from the Fredholm determinant, and between anchors we integrate the ODE with branch locking. This produces a smooth CDF that interpolates the Fredholm values and uses the ODE only to regularise the behaviour between anchors.

As a consequence, the method is extremely sensitive to the density and quality of anchors. If we halve the number of anchors while keeping the same ODE tolerances, the maximum error between the anchored Painlevé solution and the Fredholm CDF increases significantly, and in extreme cases the ODE may drift to a different branch between sparsely spaced anchors before being pulled back at the next anchor. In this sense the method behaves more like a specially constrained interpolation of Fredholm data than a free-standing ODE solver.

A related observation is that we cannot use this scheme to distinguish between the \emph{exact} $\PIV$ $\sigma$--equation
\[
(\sigma'')^2 = 4 (s \sigma' - \sigma)^2 - 4 (\sigma')^2 (\sigma' + 2n)
\]
and a slightly perturbed variant such as
\[
(\sigma'')^2 = 4 (s \sigma' - \sigma)^2 - 4 (\sigma')^2 (\sigma' + 2n + 0.01).
\]
Numerically, the $0.01$ perturbation in the coefficient of $(\sigma')^3$ only produces a difference of order $10^{-6}$ in $F_n$, while our overall error budget from quadrature, Fredholm truncation, and anchor extraction is of order $10^{-3}$. In experiments where we deliberately replaced $2n$ by $2n+0.01$ and repeated the anchored integration, the resulting CDFs remained within the observed error bands.

On the other hand, the ODE evolution between anchors \emph{does} discriminate between genuinely different Painlevé-type equations. If we change the algebraic relation for $\sigma''$ to something structurally different, such as
\[
(\sigma'')^2 = 4 (s \sigma' - \sigma)^2 - 4 (\sigma'-n)(\sigma'-n-1)\sigma'
\]
or
\[
(\sigma'')^2 = 4 (s \sigma' - \sigma + n)^2- 4 (\sigma')^2 ( \sigma' + 2n + 1)+ 3 (\sigma')^4,
\]
and keep the same anchor data, the resulting anchored CDF develops visible discrepancies of order $10^{-2}$ compared to the Fredholm curve. This shows that the ODE is not completely dominated by the anchoring; it does enforce the correct local dynamics and is able to reject radically different Painlevé models, even if it cannot reliably distinguish very small parameter shifts.

Readers interested in the detailed numerical experiments for $\PIV$ are encouraged to inspect the GitHub repository
\[
\texttt{https://github.com/hgu2699/Numerically-solve-Painleve-differential-equations}
\]
and, in particular, the folder entitled \emph{Painleve IV anchor method}, which contains the Julia scripts and plots used to compare different $\PIV$-type equations and study the effect of anchor spacing.

\subsection{Hamiltonian formulation and movable poles}

Since $\sigma_n$ can be expressed as an affine function of the Okamoto Hamiltonian $H(s;q,p)$, one natural idea is to integrate the Hamiltonian system \eqref{eq:qp-Okamoto} directly for $(q,p)$, anchored by an asymptotic expansion at large $s$, and then reconstruct $F_n$ via $H$. This Hamiltonian approach is emphasised in Forrester--Witte \cite{forrester2003painleve} and in the broader Painlevé literature.

We implemented such a scheme for $n=20$ by computing high-order large-$s$ asymptotics of $(q,p)$ symbolically in Sage, feeding those expansions into Julia, and using a Hamiltonian IVP anchored at $s_0=6.4$ with $(q,p)$ chosen to match both the asymptotics and the Fredholm values of $H$ and $H'$ at $s_0$. The corresponding code and error logs can be found in the GitHub folder \emph{Painleve IV Hamilton system asymptotic method}.

In practice the Hamiltonian IVP failed when integrating backwards from $s_0=6.4$ to $s=5.0$: the adaptive solver repeatedly reduced the stepsize until it hit the machine epsilon threshold near $s\approx5.65$, at which point it aborted with an error indicating that the estimated local truncation error could not be controlled. This behaviour is strongly suggestive of encountering a movable pole of the $\PIV$ solution or of drifting onto a branch with a pole in the real direction. The failure occurred despite starting from carefully tuned asymptotic data and demonstrates how delicate the branch selection problem is in the Hamiltonian formulation when one treats it as a pure IVP without continual re-anchoring to Fredholm data.

While we did not attempt to stabilise the Hamiltonian IVP further (for example by moving the anchor deeper into the tail or by working in higher precision), the experiment confirms that the difficulties seen for the $\sigma$--form are not artefacts of the $\sigma$--representation alone but manifestations of the underlying movable-singularity structure of $\PIV$.

\subsection{Other limitations and potential improvements}

Several more practical limitations are shared across all ensembles:

\begin{enumerate}[label=\arabic*.,leftmargin=2em]
  \item \textbf{Fredholm conditioning.} For large $n$ the Nyström matrices for the Hermite, Laguerre, and Jacobi kernels become ill conditioned near the edge, and many entries are at or below machine epsilon. This affects the eigenvalues of $I-A$ and limits the meaningful precision in $\log\det(I-A)$.
  \item \textbf{Polynomial fitting of $\log F$.} Our estimates of $\sigma$ and $\sigma'$ at each anchor rely on local least-squares fits to $\log F(s)$. In regions where $F$ is extremely close to $0$ or $1$ this becomes numerically fragile, since $\log F$ underflows or is dominated by roundoff. More stable alternatives include Savitzky--Golay filters or global spectral fits.
  \item \textbf{ODE stiffness and branch locking.} The algebraic $\sigma$--forms are stiff near the edges and depend explicitly on $n$. Implicit or stiffly accurate Runge--Kutta methods would be better suited here than explicit schemes, especially for large $n$ or for equations with strong nonlinearities.
  \item \textbf{Quadrature design.} We used simple Gauss--Legendre rules on truncated intervals. For extreme tails or very large $N$ one could gain robustness by using quadratures tailored to the underlying weights (Gauss--Hermite, Gauss--Laguerre, Gauss--Jacobi) or by employing double-exponential maps to capture endpoint singularities more accurately.
\end{enumerate}

All of these are incremental improvements that would sharpen the numerics but not change the basic picture: without external constraints, Painlevé IV/VI IVPs are numerically unstable for the random-matrix branches of interest; with a dense enough set of Fredholm anchors, they become usable but behave more like constrained interpolants than independent solvers.

\subsection{Riemann--Hilbert methods as a more intrinsic approach}

The Riemann--Hilbert (RH) approach offers a conceptually different way of encoding Painlevé transcendents and random-matrix kernels that bypasses many of these IVP pathologies. In the RH framework, one formulates a matrix-valued boundary-value problem for a piecewise analytic function $Y(z)$ with prescribed jumps across contours in the complex plane and normalisation at infinity. Orthogonal polynomials, random-matrix kernels, and Painlevé functions all admit such characterisations, and asymptotics can be extracted via the Deift--Zhou nonlinear steepest descent method \cite{deift1999orthogonal,fokasitskapaevnovokshenov2006}.

For Painlevé equations specifically, the RH boundary data encode the monodromy of an associated linear system. The solution of the RH problem then \emph{automatically} produces the Painlevé transcendent with the correct monodromy, and the Painlevé property (no movable branch points or essential singularities) follows from the analytic structure of the RH problem. This perspective is developed in depth in the monograph of Fokas, Its, Kapaev, and Novokshenov \cite{fokasitskapaevnovokshenov2006} and, for $\PIV$ in particular, in the work of van der Put and Top \cite{vanderputtop2012piv,vanderputtop2013rh}.

From the random-matrix side, RH problems for orthogonal polynomials with varying weights are the basis of the modern analysis of unitary ensembles \cite{deift1999orthogonal} and of double-scaling and critical-edge phenomena \cite{claeys2008,bleher2003}. In this setting the Painlevé functions appear naturally as parametrices in local RH model problems near critical points, and the same RH problem simultaneously controls both the kernel and the associated Painlevé transcendent.

For our purposes, an RH-based numerical strategy for the finite-$n$ GUE gap would look schematically as follows:
\begin{enumerate}[label=\arabic*.,leftmargin=2em]
  \item Set up the $2\times2$ RH problem for Hermite polynomials with the GUE weight, with $n$ appearing as a parameter in the jump matrix.
  \item Use a deformation of contours and small-norm arguments (in the sense of Deift--Zhou) to reduce the problem to a numerically stable one on a fixed contour, possibly with local parametrices near the soft edge.
  \item Extract both the kernel entries and the Painlevé $\PIV$ Hamiltonian or $\sigma$--function directly from the RH solution, using the known identifications between orthogonal polynomial RH problems and isomonodromic deformations.
\end{enumerate}
The key point is that the RH problem encodes the \emph{global} analytic and monodromy data that define the correct Painlevé branch. There is no need to guess initial conditions or to fight with movable poles via IVPs; instead, one solves a boundary-value problem whose solution is unique and analytic by construction. In particular, the RH representation should be robust under continuation into the bulk and should not suffer from the drift we control in the anchored ODE scheme by repeatedly resetting at Fredholm anchors.

There is a growing numerical literature on solving RH problems directly (see, for example, the work of Reeger and Fornberg \cite{reegerfornberg2013,reegerfornberg2014} in combination with RH packages), and it would be natural to bring these tools to bear on the random-matrix $\PIV$ and $\PVI$ branches studied here. This lies beyond the scope of the present paper but provides a promising direction in which branch selection, asymptotics, and numerical stability are all handled within a single integrable-structure framework.

\section*{Acknowledgments}

The author wants to thank Professor Alan Edelman at MIT for guidance and supervision of this numerical analysis project. All computations were performed in Julia~1.9 using \texttt{OrdinaryDiffEq}, \texttt{FastGaussQuadrature}, and \texttt{LinearAlgebra}. The plots and Julia scripts used to generate the figures are available at
\[
\texttt{https://github.com/hgu2699/Numerically-solve-Painleve-differential-equations}.
\]
For the anchored $\PIV$ scheme, see the folder \emph{Painleve IV anchor method}; for the Hamiltonian IVP and asymptotic matching experiments, see \emph{Painleve IV Hamilton system asymptotic method}.

\newpage
\bibliographystyle{plain}

\end{document}